\documentclass[12pt]{article}%,a4paper
\usepackage{graphicx,bm,epsf,amsmath,amsthm,amscd,amsfonts,amssymb,amsbsy}%,color
\usepackage[hyphens,spaces,obeyspaces]{url}%must be first!
\usepackage[colorlinks,breaklinks,allcolors=blue]{hyperref}%[dvips,
\usepackage[hdvips]{hypdvips} %!important:! - makes broken links active
\usepackage{doi}
\newcommand{\R}{\mbox{Re\,}}
\newcommand{\be}{\begin{equation}}
\newcommand{\ee}{\end{equation}}
\newcommand{\bea}{\begin{eqnarray}}
\newcommand{\eea}{\end{eqnarray}}
\newcommand{\nn}{\nonumber}
\newcommand{\noi}{\noindent}
\newcommand{\bc}{\begin{center}}
\newcommand{\ec}{\end{center}}
\newcommand {\mm}[1]{\quad\mbox{#1}\quad}
\newcommand {\MM}[1]{\qquad\mbox{#1}\qquad}

\begin{document}
\title{\bf A hypergeometric transformation formula of Ramanujan,\\
its validity range and implications}
\author{{\bf M. A. Shpot}\\[3mm]%$^{\bm a}$
{\em Yukhnovskii Institute for Condensed Matter Physics}\\
{\em of the National Academy of Sciences of Ukraine, 79011 Lviv, Ukraine}\\[2mm]%${}^a\!$
{\small\bf e-mail:} shpot.mykola@gmail.com\\[1mm]
{\small\textbf{ORCID ID:}} 0000-0001-9098-9745}

\date{\today}
\maketitle

\begin{abstract}
\noindent
We extend the validity range of a hypergeometric transformation formula of Ramanujan proved by
Berndt, Bhargava and Garvan [{\em Trans. Amer. Math. Soc.} {\bf 347} 4163 (1995)]
and study its implications.
Relations to special values of complete elliptic integrals of the first kind in
singular value theory are established.
Consequently, we derive several closed-form evaluations of hypergeometric functions $_2F_1$ with different sets of parameters and arguments.
Connections with other hypergeometric transformations and some recent results are discussed.
\end{abstract}

\medskip


\section{Introduction}

At the top of p. 254 of his handwritten notes \cite[p. 258]{Rnotes}, Ramanujan records a hypergeometric transformation relating two Gauss hypergeometric functions with \emph{different} sets of parameters
and arguments, $_2F_1(\frac13,\frac23;1;\alpha)$ and $_2F_1(\frac12,\frac12;1;\beta)$.
On that page, we find neither the proof nor indication of the validity range of the proposed relation.

In 1995, Berndt, Bhargava and Garvan \cite{BBG95} published the proof of Ramanujan's statement and formulated it as Theorem 5.6 (p. 258) \cite[p. 4184]{BBG95}, \cite[p. 112]{Berndt5}, \cite[(2.21)]{Garvan95}\footnote{This last reference presents a concise list of related hypergeometric identities on pp. 522--524.}.

It is perhaps the most famous relation among the large family of
hypergeometric transformations due to Ramanujan that follow from his
theories of elliptic functions to alternative bases (see \cite{BBG95},
\cite[Chap. 33]{Berndt5}, and \cite{Garvan95}).
We reproduce it in the following form, which is a fusion of \cite[p. 258]{Rnotes}
and \cite[p. 4184]{BBG95}:
\begin{flalign}\label{RAB}&
\mbox{{\bf Theorem 5.6} (p. 258). \emph{If}}\qquad \alpha(p):=\frac{p^3(2+p)}{1+2p}\MM{\emph{and}}
\beta(p):=\frac{27p^2(1+p)^2}{4(1+p+p^2)^3}\,,
&\end{flalign}
\emph{then, for} $0\le p<1$,
\be\label{RRR}
_2F_1\Big(\frac13,\frac23;1;\beta(p)\Big)=
\gamma(p)\,_2F_1\Big(\frac12,\frac12;1;\alpha(p)\Big)
\MM{\emph{where}}
\gamma(p):=\frac{1+p+p^2}{\sqrt{1+2p}}\,.
\ee
For brevity, we shall refer to this as the RBBG relation.

A corollary of the above theorem is the complementary identity \cite[(5.21)]{BBG95},
\be\label{Cor}
\sqrt3\,_2F_1\Big(\frac13,\frac23;1;1-\beta(p)\Big)=
\gamma(p)\,_2F_1\Big(\frac12,\frac12;1;1-\alpha(p)\Big).
\ee
\medskip

The function $_2F_1(\frac13,\frac23;1;z)$ is a central ingredient
of Ramanujan's theories of elliptic functions to alternative bases
from his notebooks \cite{Berndt5} and attracted considerable attention in subsequent works.
In particular, the RBBG formula \eqref{RRR} has been proven in different ways
in \cite[p. 201]{Chan98}, \cite[p. 357]{CZ19}, \cite{Rob21} between 1998 and 2021.

A companion version of the RBBG identity
\be\nn
_2F_1\Big(\frac13,\frac23;1;\tilde\beta(p)\Big)=
\tilde\gamma(p)\,_2F_1\Big(\frac12,\frac12;1;\alpha(p)\Big)
\;\;\mbox{with}\;
\tilde\beta(p)=\frac{27p(1+p)^4}{2(1+4p+p^2)^3}\,,\;
\tilde\gamma(p)=\frac{1+4p+p^2}{\sqrt{1+2p}}
\ee
can be found in recent publications \cite[p. 5]{RingZ} and \cite[p. 35]{Ring23}.
\medskip

In what follows, we shall show that the RBBG relation \eqref{RRR}
holds over the extended real interval $-\frac12<p<1$.
Throughout this range, $\beta(p)\in[0,1[$ and $\alpha(p)\in]-\infty,1[$.
At the endpoints $p=-\frac12$ and $p=1$, both $_2F_1$ functions in \eqref{RRR} diverge
(see Sec.~\ref{SCR} for further details).

\section{Selected properties of the Gauss $_2F_1$ function}\label{CGF}

For convenience, before proceeding to the next section, let us recall some basic facts
concerning the convergence properties of the Gauss hypergeometric series
\be\label{HYS}
_2F_1(a,b;c;z)=\sum_{n\ge0}\,\frac{(a)_n(b)_n}{(c)_n}\,\frac{z^n}{n!}\,,
\qquad c\ne0,-1,-2\,,\ldots
\ee
where the Pochhammer symbol $(\lambda)_n$ \cite{Pochhammer} is given by
\[
(\lambda)_n=\dfrac{\Gamma(\lambda+n)}{\Gamma(\lambda)}=
\begin{cases}\lambda(\lambda+1)\ldots (\lambda+n-1)&n\in\mathbb N;\;\\
             1                   &n=0\,;
\end{cases}\;\;\lambda\in\mathbb C\,.
\]

1$^{\bm\circ}$\, For any $z$ inside the unit circle $|z|<1$, the hypergeometric series \eqref{HYS} converges absolutely provided that
all its parameters $a$, $b$ and $c$ differ from zero and negative integers.%
\footnote{\label{fnp} If either of the numerator parameters $a$ or $b$ is zero or a negative integer, the series \eqref{HYS} terminates and the question of convergence does not apply.}

\medskip
2$^{\bm\circ}$\, On the unit circle $|z|=1$, the series \eqref{HYS}
\medskip

\noi($i$)\;\;converges absolutely if the real part of the parametric excess $s:=c-a-b$ is positive: $\R s>0$.

\medskip
\emph{In this case, if $z=1$, the sum of the hypergeometric series is given by
Gauss's theorem}
\be\label{GTE}
_2F_1(a,b;c;1)=\frac{\Gamma(c)\Gamma(s)}{\Gamma(c-a)\Gamma(c-b)}\,,
\qquad c\ne0,-1,-2\,,\ldots\,,\quad s=c-a-b,\quad\R s>0.
\ee

\noi($ii$)\;\;converges conditionally if $-1<\R s\le0$ and $z\ne1$.

\medskip
\emph{This case includes the point $z=-1$ where the sum of the series can be determined by Kummer's theorem} \cite[p. 9, (1)]{Bailey}, \cite[7.3.6.2]{PBM3}
\be\label{KTE}
_2F_1(a,b;a-b+1;-1)=\frac{\Gamma(a-b+1)\Gamma(1+\frac12\,a)}
{\Gamma(1+a)\Gamma(1+\frac12\,a-b)}=
2^{-a}\sqrt\pi\,\frac{\Gamma(a-b+1)}{\Gamma\big((1+a)/2\big)\Gamma(1+\frac12\,a-b)}\,.
\ee

\medskip
\noi($iii$)\;\;diverges for all $z$ with $|z|=1$ if $\R s\le-1$.

\medskip
3$^{\bm\circ}$\, If $|z|>1$, the hypergeometric series \eqref{HYS} diverges provided that it
does not terminate ($a,\;b\ne0,-1,-2\,,\ldots$).$^{\ref{fnp}}$

\medskip\noi
For further details and proofs, we refer to standard textbooks
\cite{bateman,Slater,Rainville,AAR,NIST}.

Beyond the convergence region of the series \eqref{HYS}, the Gauss hypergeometric function $_2F_1(a,b;c;z)$ is defined by analytic continuation.

One of the simplest possibilities of the analytic continuation is provided by the
linear Pfaff--Kummer transformation
\cite[p. 76]{bateman}, \cite[p.~60]{Rainville}, \cite[p.~68]{AAR}
\be\label{Lin}
_2F_1(a,b;c;z)=(1-z)^{-a}\,{}_2F_1\Big(a,c-b;c;\frac z{z-1}\Big).
\ee
This result, initially valid for both $|z|<1$ and $|z/(z-1)|<1$, can be used
for analytic continuation of $_2F_1(a,b;c,z)$ into the half-plane $\R z<\frac12$
since the right-hand side of \eqref{Lin} converges in this region.
For real $z$, which is of primary interest in this paper, the transformation
\eqref{Lin} maps the interval $z\in[0,\frac12]$ onto the segment $z\in[-1,0]$
and the interval $z\in[\frac12,1]$ is mapped onto the ray $z\in]-\infty,-1]$.
Thus, for $_2F_1$ functions with large negative arguments $z\in]-\infty,-1]$,
the analytic continuation \eqref{Lin} yields real values expressed in terms of those with arguments in the interval $[\frac12,1]$.

\section{The validity range of the RBBG equation}\label{VRR}

Having proven Ramanujan's relation \eqref{RRR}, Berndt, Bhargava and Garvan
\cite[p. 4185]{BBG95} offer the following:

\noi
\emph{... we show that our proof above of \eqref{RRR} is valid for $0<p<1$.}

\emph{Observe that}
\be\label{ders}
\frac{d\alpha}{dp}=\frac{6p^2(1+p)^2}{(1+2p)^2}\MM{and}
\frac{d\beta}{dp}=\frac{27p(1-p^2)(1+2p)(2+p)}{4(1+p+p^2)^4}\,.
\ee
\emph{Hence, $\alpha(p)$ and $\beta(p)$ are monotonically increasing on $(0,1)$.
Since $\alpha(0)=0=\beta(0)$ and $\alpha(1)=1=\beta(1)$, it follows that
Theorem 5.6 is valid for $0<p<1$.}

This observation implies that if $0<p<1$, both arguments of $_2F_1$
functions in \eqref{RRR},
$\alpha(p)$ and $\beta(p)$, lie within the interval $(0,1)$.
Hence, by property 1$^{\bm\circ}$ of the preceding section, the involved hypergeometric functions do converge here and thus a necessary condition for the
validity of \eqref{RRR} on $(0,1)$ is fulfilled.

While Berndt, Bhargava and Garvan establish that relation \eqref{RRR} holds for $0<p<1$, this does not restrict its validity strictly to this interval.
We shall show that both the necessary convergence condition and the full validity range of the RBBG relation extend to a larger interval of $p$.

\subsection{Convergence of hypergeometric functions in \eqref{RRR}}\label{SCR}

Both functions $_2F_1(a,b;c;z)$ in \eqref{RRR} are zero-balanced: their parametric excess $s=c-a-b$ vanishes. Thus, by properties  1$^{\bm\circ}$
and 2$^{\bm\circ}$($ii$) of Sec.~\ref{CGF}, they converge for arguments $z\in[-1,1[$.

\begin{figure}[htb]
\bc
\includegraphics[width=8cm,height=6cm]{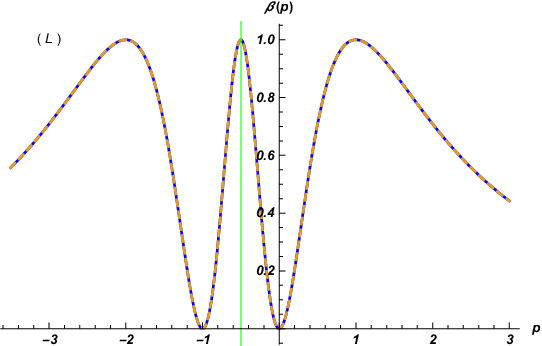}
\includegraphics[width=8cm,height=6cm]{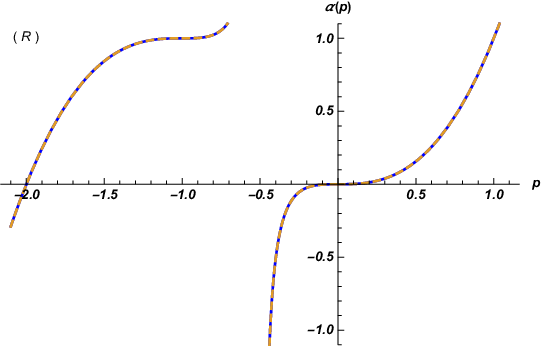}%\\\hspace{5mm}(L)\hspace{80mm}(R)
\ec
\vspace{-5mm}
\caption{(L) Coinciding functions $\beta(p)$ (blue dots) and $\beta(p^{-1})$ (dashed orange); the symmetry axis $p=-\frac12$ is shown as a green vertical line.
(R) Coinciding functions $\alpha(p)$ (blue dots) and $1/\alpha(p^{-1})$ (dashed orange).}
\label{PB}
\end{figure}

The argument $\beta(p)$ is plotted in Fig.~\ref{PB}(L). The function $\beta(p)$ is
non-negative and bounded within the interval $[0,1]$. Consequently, the series
$_2F_1(\frac13,\frac23;1;\beta(p))$ converges for all $p$ except those for which
$\beta(p)=1$, namely at $p=-2\,,-\frac12\,,1$.

The function $\beta(p)$ satisfies the following symmetry properties:
\vspace{-2mm}
\begin{itemize}\itemsep-1mm
\item It is invariant under inversion $p\leftrightarrow p^{-1}$, that is,
$\beta(p)=\beta(1/p)$;
\item It is symmetric with respect to the vertical line $p=-\frac12$ (green line
in Fig.~\ref{PB} (L)).
\end{itemize}

The behavior of the function $\alpha(p)$ is less straightforward --- see Fig.~\ref{PB}(R). It diverges as $|p\,|\to\infty$ and $p\to-\frac12$ exhibiting an  infinite discontinuity at $p=-\frac12$ .
In addition, it satisfies the inversion relation $\alpha(p)=1/\alpha(p^{-1})$.

While $p\in]0,1[$ has been considered in \cite{BBG95}, only the corresponding
positive portion of the argument $\alpha(p)\in]0,1[$ was taken into account.
However, the function $_2F_1(\frac12,\frac12;1;\alpha(p))$ still converges for
\emph{negative}\footnote{In the special case $p=0$, relation \eqref{RRR} holds trivially.}
$\alpha(p)$, corresponding to $p\lesssim0$, provided that $\alpha(p)\ge-1$.
Solving the equation $\alpha(p)=-1$ using Mathematica \cite{Math12}, we obtain the special value of $p$,
\be\label{pesc}
p^*=\frac12\left(-1+\sqrt{2\sqrt3}-\sqrt3\right)\simeq-0.435,
\ee
below which $|\alpha(p)|>1$ and the function $_2F_1(\frac12,\frac12;1;\alpha(p))$
ceases to converge.
Thus, $\alpha(p)\in[-1,1[$ for $p\in[p^*,1[$, and the function
$_2F_1(\frac12,\frac12;1;\alpha(p))$ converges for any $p$ in this interval.

\medskip
Thus, we obtain a simple extension of the necessary condition for the validity of
relation \eqref{RRR}:

\noi
\emph{For the RBBG relation to hold, it is necessary that the parameter $p$ must belong
to the interval $p^*\le p<1$, since both hypergeometric functions involved in \eqref{RRR} converge in this interval.}
\medskip

\begin{figure}[htb]
\bc
\includegraphics[width=8cm,height=6cm]{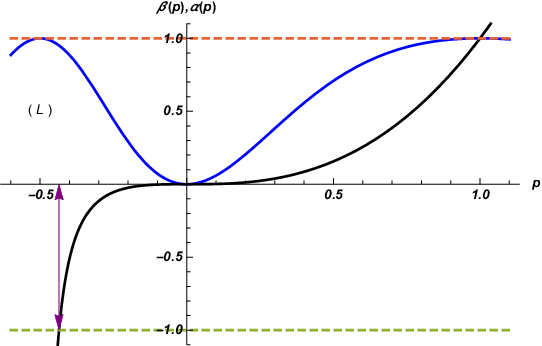}
\includegraphics[width=8cm,height=6cm]{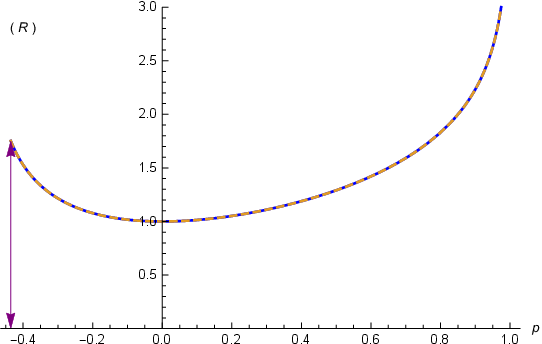}%\\\hspace{5mm}(L)\hspace{80mm}(R)
\ec
\vspace{-5mm}
\caption{(L) Functions $\beta(p)$ (blue) and $\alpha(p)$ (black) in the interval
$p^*\le p<1$ with $p^*$ given in \eqref{pesc}, and slightly beyond.
Dashed lines at $\pm1$ indicate the convergence boundaries for the $_2F_1$ functions in \eqref{RRR}.
The series converge along the green boundary line ($\alpha=-1$), whereas they diverge along the red line ($\alpha,\beta=1$).
The purple arrows indicate  $(i)$ the point $p^*$ on the $p$-axis and
$(ii)$ the boundary point where $\alpha(p)$ leaves the convergence domain of $_2F_1(\frac12,\frac12;1;\alpha(p))$.\\
(R) Coinciding left-hand (blue dots) and right-hand (dashed orange) sides of \eqref{RRR}. The purple arrows again indicate $p^*$ on the abscissa and the corresponding value of plotted functions.
}
\label{P2}
\end{figure}

Figure~\ref{P2}(L) shows the arguments $\alpha(p)$ and $\beta(p)$ taking values in $[-1,1]$ over the interval $p^*\le p<1$ and its immediate vicinity.
Throughout this interval, the RBBG relation \eqref{RRR} is well-defined, as both arguments of the respective $_2F_1$ series lie within their convergence domains.
A numerical verification confirms that \eqref{RRR} is valid for $p^*\le p<1$ with $p^*$ given in \eqref{pesc}, as shown in Fig.~\ref{P2}(R).

Below, we explicitly calculate the values of both $_2F_1$ functions in \eqref{RRR} at $p=p^*$.

\subsection{Further extension of the validity range of the RBBG formula}

\begin{figure}[htb]
\bc
\includegraphics[width=8cm,height=6cm]{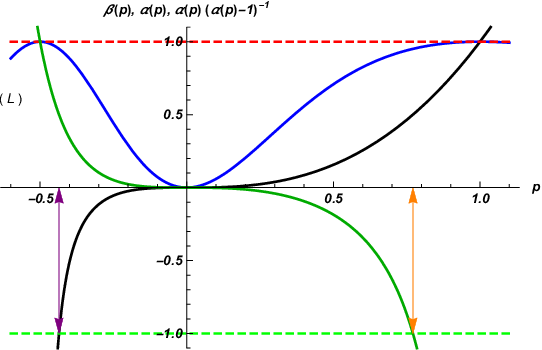}
\includegraphics[width=8cm,height=6cm]{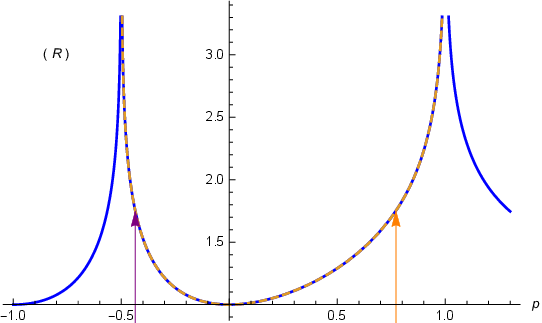}%\\\hspace{5mm}(L)\hspace{80mm}(R)
\ec
\vspace{-5mm}
\caption{(L) Functions $\beta(p)$ (blue), $\alpha(p)$ (black), and $\alpha_\ell(p)=\alpha(p)/(\alpha(p)-1)$ (dark green) over the interval of validity $p\in]-\frac12,1[$ of relation \eqref{RRR}. As before, the purple arrows indicate $p=p^*$ and the point $\alpha(p^*)=-1$ where the argument $\alpha(p)$ leaves the convergence domain. The orange arrows show the value $p^*_\ell$ and the corresponding boundary point for $\alpha_\ell(p)$.
(R) Coincidence of the functions $_2F_1(\frac13,\frac23;1;\beta(p))$ (blue) and $\gamma(p)\,_2F_1(\frac12,\frac12;1;\alpha(p))$ (dashed orange) over $p\in]-\frac12,1[$. The purple and orange arrows indicate their common value $1.7514579\ldots$ at $p=p^*$ and $p=p^*_\ell$, as defined in \eqref{pesc} and \eqref{pell}, respectively.}
\label{P3}
\end{figure}

When plotting the function $_2F_1(\frac12,\frac12;1;\alpha(p))$,
Mathematica automatically employs its analytic continuation yielding real values appropriate for visualization.
This provides clear numerical evidence that the validity range of relation \eqref{RRR} extends further to $-\frac12<p<1$, as shown in Fig.~\ref{P3}(R).

Explicitly, such analytic continuation is achieved by applying the Pfaff-Kummer linear transformation \eqref{Lin} to the right-hand side of the RBBG formula \eqref{RRR}. This results in the following modification of the original equation:
\be\label{RRL}
_2F_1\Big(\frac13,\frac23;1;\beta(p)\Big)=
\gamma_\ell(p)\,_2F_1\Big(\frac12,\frac12;1;\alpha_\ell(p)\Big)
\ee
where
\be\label{GAL}
\gamma_\ell(p):=\frac{1+p+p^2}{(1+p)\sqrt{1-p^2}}
\MM{and}
\alpha_\ell(p):=\frac{-p^3(2+p)}{(1-p^2)(1+p)^2}\,.
\ee

As described at the end of Sec.~\ref{CGF}, the infinite segment of the function
$\alpha(p)$ below the boundary $\alpha(p)=-1$ (see Fig.~\ref{P3}(L)) maps onto
the interval of the transformed argument $\frac12<\alpha_\ell(p)<1$.
This occurs in the interval $p\in]-\frac12,p^*\,[$, complementary to the one
considered in Figure~\ref{P2}.

The hypergeometric series $_2F_1(\frac12,\frac12;1;\alpha_\ell(p))$ converges for all $p\in]-\frac12,p^*_\ell]$ and ceases to converge when $p$ exceeds the special value $p^*_\ell$,
\be\label{pell}
p^*_\ell=\frac12\left(-1+\sqrt{2\sqrt3+3}\right)=
\frac12\sqrt{\frac{\sqrt3}2}\left(1+\sqrt3-\sqrt{\frac2{\sqrt3}}\right)\simeq0.771,
\ee
at which the transformed argument $\alpha_\ell(p)$ equals $-1$.
The value $p^*_\ell$ is shown by the upward orange arrow in Fig.~\ref{P3}(L),
while the downward orange arrow indicates the point where $\alpha_\ell(p)$ leaves the convergence domain between the red and green dashed lines.

Although the analytic expressions on the right-hand sides of \eqref{RRR} and \eqref{RRL} differ, the plot for the modified RBBG equation \eqref{RRL} is identical to that shown in Fig.~\ref{P3}(R).
Mathematica makes no distinction between the functions on the right-hand sides of equations \eqref{RRR} and \eqref{RRL}.

In summary, we formulate the following proposition, which is an extended version of Theorem 5.6 by Berndt, Bhargava, and Garvan:

\vspace{1mm}
\noi
\emph{If}
\be\nn
\beta(p)=\frac{27p^2(1+p)^2}{4(1+p+p^2)^3}\,,\qquad
\alpha(p)=\frac{p^3(2+p)}{1+2p}\,,\qquad
\alpha_\ell(p)=\frac{-p^3(2+p)}{(1-p^2)(1+p)^2}\,,
\ee
\be\nn
\gamma(p)=\frac{1+p+p^2}{\sqrt{1+2p}}
\MM{\emph{and}}
\gamma_\ell(p)=\frac{1+p+p^2}{(1+p)\sqrt{1-p^2}}\,,
\ee
\emph{then, for $-\frac12<p<1$,}
\be\label{BRR1}
_2F_1\Big(\frac13,\frac23;1;\beta(p)\Big)=
\left\{\begin{array}{l}\displaystyle{
\gamma_\ell(p)\,_2F_1\Big(\frac12,\frac12;1;\alpha_\ell(p)\Big)
\;\mm{\emph{for}}-\frac12<p\le p^*_\ell,
}
\medskip\\
\displaystyle{
\gamma(p)\,_2F_1\Big(\frac12,\frac12;1;\alpha(p)\Big)
\MM{\emph{for}}p^*\le p<1,
}
\end{array}\right.
\ee
\emph{where $p^*$ and $p^*_\ell$ satisfy $-\frac12<p^*<0<p^*_\ell<1$ and are given by}
\be\nn
p^*=-\frac12\left(1+\sqrt3-\sqrt{2\sqrt3}\right)
\MM{\emph{and}}
p^*_\ell=\frac12\left(\sqrt{3+2\sqrt3}-1\right).
\ee
\emph{In each case, under the indicated restrictions on $p$, the hypergeometric series in \eqref{BRR1} converge.
In the overlapping interval $p^*\le p\le p^*_\ell$, both series
$_2F_1(\frac12,\frac12;1;\alpha(p))$ and
$_2F_1(\frac12,\frac12;1;\alpha_\ell(p))$ converge simultaneously.}
\vspace{1mm}

In support of our statement, in the following sections we present analytical calculations for certain special values of $p$, such as the boundary points $p^*$ and $p^*_\ell$ from \eqref{pesc} and \eqref{pell}, and discuss their implications.

\section{Evaluations at the boundary points}

The following evaluations provide additional evidence that the RBBG formula holds throughout the extended validity domain established in Sec.~\ref{VRR}.

As discussed in Sec.~\ref{SCR}, the point $p^*<0$ is the marginal value of the variable $p$ at which the hypergeometric series $_2F_1(\frac12,\frac12;1;\alpha(p))$
on the right-hand side of the RBBG formula \eqref{RRR} still converges.
At $p=p^*$ given in \eqref{pesc}, the argument of this $_2F_1$ function is $\alpha(p^*)=-1$.
Thus, its sum is given by Kummer's theorem \eqref{KTE}, yielding
\be
_2F_1\Big(\frac12,\frac12;1;-1\Big)=
\frac2{\sqrt\pi}\,\frac{\Gamma(\frac54)}{\Gamma(\frac34)}=
\frac{\Gamma^2(\frac14)}{(2\pi)^{3/2}}\,.
\ee
Further, using the original definition of the function $\beta(p)$ given in \eqref{RAB}, we compute
\be\nn
\beta(p^*)=\beta(p^*_\ell)=1-\beta(p_9)=\frac{3\sqrt3}4\,\big(\sqrt3-1\big)\simeq0.951
\MM{and}
\gamma(p^*)=\Big(\frac{3\sqrt3}2\Big)^{1/4}\sqrt{\sqrt3+1}\;.
\ee
The values $p^*$ and $p^*_\ell$ are given in \eqref{pesc} and \eqref{pell},
while $p_9$ is defined as
\be
p_9=\frac12\left(\sqrt{\frac{3\sqrt3}2}\left(\sqrt3-1\right)-1\right).
\ee
In a different form, the same value of $p$ will appear shortly in the
next section, see \eqref{P9}.
The explicit value of $\beta(p_9)$ is given in \eqref{IIN}, and
in Sec.~\ref{BFS} it will be denoted for brevity as $y_1$, see \eqref{Y0}.

Substituting these results into the RBBG relation yields
\be\label{comm}
_2F_1\Big(\frac13,\frac23;1;1-y_1\Big)=
\Big(\frac{3\sqrt3}2\,\Big)^{1/4}\,\frac{\sqrt{\sqrt3+1}}{(2\pi)^{3/2}}\;
\Gamma^2\Big(\frac14\Big),
\ee
where
\be
1-y_1=1-\Big(\frac{\sqrt3-1}2\Big)^3=\frac{3\sqrt3}4\,\big(\sqrt3-1\big).
\ee

The result \eqref{comm} is related to the Berndt--Chan--Ramanujan evaluation expressed in two equivalent forms in \eqref{BR1} and \eqref{BR3}, via Ramanujan's cubic transformation
\cite[Chap.~35, (4.21)]{Berndt5}
\be\label{RCT}
_2F_1\Big(\frac13,\frac23;1;1-\Big(\frac{1-x}{1+2x}\Big)^3\Big)=
(1+2x)\,_2F_1\Big(\frac13,\frac23;1;x^3\Big).
\ee
Dividing \eqref{comm} by \eqref{BR3}, we reproduce the result
\be
\frac{\displaystyle _2F_1\Big(\frac13,\frac23;1;\frac{3\sqrt3}4\,\big(\sqrt3-1\big)\Big)}
{\displaystyle_2F_1\Big(\frac13,\frac23;1;\frac{3\sqrt3-5}4\Big)}=\sqrt3
\ee
from \cite[p.~292]{BC95} and \cite[Ch.~35, (4.22)]{Berndt5}.

The above example shows how the RBBG relation works for a negative value of $p$, included into its extended validity region in Sec.~\ref{SCR}.
We have also seen that the evaluation at the boundary point $p=p^*<0$ is closely related to the well-known explicit Berndt--Chan--Ramanujan results and to
the $n=9$ value of $_2F_1(\frac12,\frac12;1;x_n)$ arising in the singular-value theory of elliptic integrals of the first kind.
Both of these connections are examined in greater detail in the next section.

\section{Berndt--Chan--Ramanujan evaluations}

Let us quote two explicit determinations of $_2F_1$
appearing in the paper by Bernd and Chan \cite[p.~280]{BC95} and in Berndt's book \cite[pp.~327--8]{Berndt5} as corollaries of an entry from Ramanujan's notebooks:

\medskip
\noi
\emph{If}
\be\label{P9}
p=\frac{\sqrt{6\sqrt3-9}-1}{2}:=p_9\,,
\ee
\emph{then}
\be\label{BR2}
_2F_1\Big(\frac12,\frac12;1;\alpha(p_9)\Big)=
\frac{\sqrt\pi}{\sqrt{6\sqrt3-9}\;\Gamma^2(\frac34)}
=\frac{1+\sqrt3}{(2\sqrt3\,\pi)^{3/2}}\,\Gamma^2\Big(\frac14\Big)
\ee
\emph{and}
\be\label{BR1}
_2F_1\Big(\frac13,\frac23;1;\beta(p_9)\Big)=
\frac{\sqrt\pi}{2^{1/4}\cdot3^{1/8}\sqrt{\sqrt3-1}\;\Gamma^2(\frac34)}
\ee
\emph{where}
\be\label{IIN}
\alpha(p_9)=\left(\frac{\sqrt2-\sqrt[4]{3}}{1+\sqrt3}\right)^2:=x_9
\MM{\emph{and}}
\beta(p_9)=\frac{3\sqrt3-5}4\,,
\ee
\emph{while the functions $\alpha(p)$ and $\beta(p)$ are defined in \eqref{RAB}, as before.}
\medskip

While the middle evaluation in \eqref{BR2} was a direct consequence of an entry by Ramanujan, the result \eqref{BR1} has been derived \cite[p.~280]{BC95}
by employing Theorem \eqref{RRR}.
Thus, the functions in \eqref{BR2} and \eqref{BR1} satisfy the relation \eqref{RRR} with the coefficient
\be
\gamma(p_9)=\frac{3^{5/8}}{2^{1/4}\sqrt{\sqrt3+1}}\,.
\ee

With the identification $\alpha(p_9)=x_9$ in \eqref{IIN}, the result \eqref{BR2}
is equivalent to the evaluation of the singular value $_2F_1(\frac12,\frac12;1;x_9)$
for the $n=9$ case of the singular-value theory of elliptic integrals of the first kind \cite[p.~117]{Z77}, \cite[p.~259]{JZ91}, \cite[Appendix~A.3]{Lattice},
\cite[Sec.~4]{SR14}, see also \cite{RSS26}.
In this theory, for singular values of positive integer order $n$, one has
\be\label{CON}
\frac{_2F_1(\frac12,\frac12;1;1-x_n)}{_2F_1(\frac12,\frac12;1;x_n)}=
\frac{K(k'_n)}{K(k_n)}=\sqrt n\,,\qquad
n\in\mathbb N,
\ee
where $x_n=k_n^2$, and, as usual, the singular moduli $k_n$ and $k'_n$ of complete
elliptic integrals of the first kind $K$ are related via $k'_n=\sqrt{1-k_n^2}$.
Recall that the relation between $K$ and the Gauss hypergeometric $_2F_1$ function is
\be
K(k)=\frac\pi2\,_2F_1\Big(\frac12,\frac12;1;k^2\Big).
\ee

Expressing the original Berndt--Chan equation \eqref{BR1} in terms of $\Gamma(\frac14)$, as we have done in \eqref{BR2}, yields
\be\label{BR3}
_2F_1\Big(\frac13,\frac23;1;\frac{3\sqrt3-5}4\Big)=(2\sqrt3\,)^{-1/4}\,
\frac{\sqrt{\sqrt3+1}}{(2\pi)^{3/2}}\,\Gamma^2\Big(\frac14\Big).
\ee
This is very similar to the evaluation \cite[(48)]{RSS26} obtained via the Mingari Scarpello--Ritelli equation \cite[(12)]{SR18}, namely,
\be\label{RS3}
_2F_1\Big(\frac13,\frac23;1;\frac{-3\sqrt3-5}4\,\Big)=(2\sqrt3\,)^{-1/4}\,
\frac{\sqrt{\sqrt3-1}}{2\pi^{3/2}}\,\Gamma^2\Big(\frac14\Big).
\ee

Applying the linear Pfaff--Kummer transformation \eqref{Lin} to the $_2F_1$ function in \eqref{BR3} yields the evaluation
\be\label{B33}
_2F_1\Big(\frac13,\frac13;1;\frac{\sqrt3-2}{3\sqrt3}\,\Big)=
\frac{3^{3/8}}{2^{1/12}4\pi^{3/2}}\,(\sqrt3+1)^{1/6}\,\Gamma^2\Big(\frac14\Big).
\ee
This is a companion evaluation to that of \cite[(47)]{RSS26},
\be\label{R33}
_2F_1\Big(\frac13,\frac13;1;\frac{\sqrt3+2}{3\sqrt3}\,\Big)=
\frac{3^{3/8}}{2^{1/12}(2\pi)^{3/2}}\,(\sqrt3-1)^{1/6}\,\Gamma^2\Big(\frac14\Big).
\ee
It is instructive to take the ratios of these two pairs of functions, which yields \cite[(51)]{RSS26}
\be\label{BRR}
\frac{\displaystyle_2F_1\Big(\frac13,\frac23;1;\frac{-3\sqrt3-5}4\,\Big)}
{\displaystyle_2F_1\Big(\frac13,\frac23;1;\frac{3\sqrt3-5}4\,\Big)}=\sqrt3-1
\mm{and}
\frac{\displaystyle_2F_1\Big(\frac13,\frac13;1;\frac{\sqrt3+2}{3\sqrt3}\,\Big)}
{\displaystyle_2F_1\Big(\frac13,\frac13;1;\frac{\sqrt3-2}{3\sqrt3}\,\Big)}=
2^{1/3}\Big(\sqrt3-1\Big)^{1/3}\!.
\ee
These ratios turn out to be special cases of a broader family of one-parameter evaluations related to the work of Ebisu, Beukers and Forsg{\aa}rd, which we discuss next.

\section{Some sample evaluations of Beukers and Forsg{\aa}rd}\label{BFS}

In a recent paper \cite[Sec.~6]{BF22}, Beukers and Forsg{\aa}rd presented several "sample $\Gamma$-evaluations" of certain $_2F_1$ functions whose parameters depend on a single variable $t$, while their arguments coincide with those appearing in \eqref{BRR}. These may be viewed as generalizations of the classical results considered in the preceding section.

Two pairs of such functions corresponding to the "shift vector" $\gamma=(-2,4,2)$ appear on p.~696 of \cite{BF22}. The involved combinations of
$\Gamma$-functions are identical in the first two evaluations, as well as in the similar formula at the top of p.~697.
In terms of the general Pochhammer symbol $(\lambda)_\nu:=\Gamma(\lambda+\nu)/\Gamma(\lambda)$ (see e.g., \cite{SS17}),
these $\Gamma$-function combinations can be compactly written as
\be\label{CGA}
C_1(t)=\frac{(\frac23)_t(\frac76)_t}{(\frac34)_t(\frac{13}{12})_t}\,.
\ee

Defining the parameter $a=-2t$ and using the original arguments
\begin{align}\label{ZZ0}
&z_0=\frac{\sqrt3+2}{3\sqrt3}=\frac{(\sqrt3+1)^2}{6\sqrt3}=
\frac4{3\sqrt3}\,\cos^2\frac\pi{12}\,,
\\\label{ZZ1}
&z_1=\frac{\sqrt3-2}{3\sqrt3}=-\frac{(\sqrt3-1)^2}{6\sqrt3}=
-\frac4{3\sqrt3}\,\sin^2\frac\pi{12}
\end{align}
of the first pair in \cite[p.~696]{BF22}, we write this set of evaluations as
\begin{align}&\label{BF1}
_2F_1\Big(a,1-2a;\frac43-a;z_0\Big)=\Big({-}\frac{27z_1}{16}\Big)^{\!-\frac a2}\,
\frac{\cos\frac\pi2\left(a+\frac16\right)}{\cos\frac\pi{12}}\,\left. C_1(t)\right|_{t=-a/2},
\\&\label{BF1A}
_2F_1\Big(a,1-2a;\frac43-a;z_1\Big)=\Big(\frac{27z_0}{16}\Big)^{\!-\frac a2}\,
\left. C_1(t)\right|_{t=-a/2}\,.
\end{align}
Taking into account the second form of $z_0$ in \eqref{ZZ0}, we observe that Beukers--Forsg{\aa}rd's result \eqref{BF1A} is identical to the entry~(ix) in Ebisu's book \cite[p.~72]{Ebisu}.
Likewise, evaluation \eqref{BF1} is equivalent to entry \cite[p.~69, (ix)]{Ebisu} although the latter is expressed in a slightly different form.

The similarity of evaluations \eqref{BF1} and \eqref{BF1A} prompts to
consider their ratio, which is given by
\be\label{RF1}
R_1(a):=\frac{_2F_1\left(a,1-2a;\frac43-a;z_0\right)}
{_2F_1\left(a,1-2a;\frac43-a;z_1\right)}=
2\left(\frac{\sqrt3+1}{\sqrt2}\right)^{2a-1}\cos\frac\pi2\Big(a+\frac16\Big),
\ee
representing a cosine function modulated by an exponential factor.
In the trigonometric form, this becomes
\be\label{RF1t}
R_1(a)=4^a\Big(\cos\frac\pi{12}\Big)^{2a-1}\cos\frac\pi2\Big(a+\frac16\Big).
\ee
In the special case $a=\frac13$, the ratio $R_1(\frac13)=\left[2(\sqrt3-1)\right]^{1/3}$
reproduces the second result in \eqref{BRR}.

The second pair of $_2F_1$ functions in \cite[p.~696]{BF22}, associated with the shift vector $\gamma=(-2,4,2)$, also contains identical combinations of $\Gamma$-functions similar to $C_1(t)$ in \eqref{CGA}. The corresponding ratio is
\be\label{BF2}
R_2(a):=\frac{_2F_1\left(a,2-2a;\frac53-a;z_0\right)}
{_2F_1\left(a,2-2a;\frac53-a;z_1\right)}=
2\left(\frac{\sqrt3+1}{\sqrt2}\right)^{2a-1}\cos\frac\pi2\Big(a-\frac16\Big).
\ee
Compared with the previous case \eqref{RF1}, in the $_2F_1$ functions appearing in \eqref{BF2}, one of the numerator parameters is increased by $1$, while
the denominator parameter is increased by $\frac13$.
Remarkably, these additive parameter shifts merely result in a sign change inside the argument of the cosine function.

\medskip
We have already noticed that another $_2F_1$ function in \cite[p.~697]{BF22}, corresponding to the shift vector $\gamma=(-2,-2,2)$, contains the same combination of $\Gamma$-functions \eqref{CGA} as the pair \eqref{BF1}--\eqref{BF1A}. Thus, it can be written in the form
\be\label{FF3}
_2F_1\Big(a,a+\frac13;\frac43-a;y_1\Big)=\left(\frac{81\sqrt3}{128}\right)^{\!-\frac a2}\,
\left. C_1(t)\right|_{t=-a/2}
\ee
where we introduced the argument notation $(3\sqrt3-5)/4:=y_1$ instead of $z_0$ appearing in \cite[p.~697]{BF22}.
Note that the same argument \eqref{IIN} has been encountered before, in the classical evaluations \eqref{BR1} and \eqref{BR3}, and appeared in the first ratio of \eqref{BRR}. The result \eqref{FF3} can be found as well in \cite[p.~73,~(xi)]{Ebisu}.

\begin{figure}[htb]
\bc
\includegraphics[width=10cm,height=7cm]{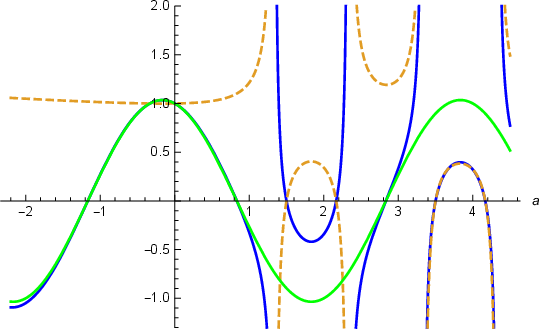}
\ec
\caption{Emergence of a cosine function (green) obtained by dividing
$_2F_1\big(a,\frac13+a;\frac43-a;y_0\big)$ (blue)
by $_2F_1\big(a,\frac13+a;\frac43-a;y_1\big)$ (dashed orange), see \eqref{BF3}.
The arguments $y_1$ and $y_0$ are given in \eqref{Y0} and \eqref{Y1}, respectively.}
\label{RRP}
\end{figure}

As a matter of fact, the Gauss function in \eqref{FF3} is merely obtained from the one in \eqref{BF1A} by applying the linear transformation \eqref{Lin}. Thus, its argument satisfies
\be\label{Y0}
y_1=\frac{3\sqrt3-5}4=\frac{z_1}{z_1-1}=\Big(\frac{\sqrt3-1}2\Big)^3=
\Big(\sqrt2\,\sin\frac\pi{12}\Big)^3
\ee
where $z_1$ is given in \eqref{ZZ1}.

Now, it is straightforward to see that a companion evaluation to \eqref{FF3} will follow from an application of the linear transformation \eqref{Lin} to the function \eqref{BF1} yielding
\be\label{LAS}
_2F_1\Big(a,a+\frac13;\frac43-a;y_0\Big)=\left(\frac{81\sqrt3}{128}\right)^{-a/2}\sqrt2\,
\big(\sqrt3-1\big)\cos\frac\pi2\Big(a+\frac16\Big)\left. C_1(t)\right|_{t=-a/2},
\ee
where
\be\label{Y1}
y_0=-\frac{3\sqrt3+5}4=\frac{z_0}{z_0-1}=-\Big(\frac{\sqrt3+1}2\Big)^3=
-\Big(\sqrt2\,\cos\frac\pi{12}\Big)^3.
\ee

Dividing equation \eqref{LAS} by \eqref{FF3} we obtain the ratio
\be\label{BF3}
R_3(a)=
\frac{_2F_1\big(a,\frac13+a;\frac43-a;y_0\big)}
{_2F_1\big(a,\frac13+a;\frac43-a;y_1\big)}=
\frac{\cos\frac\pi2(a+\frac16)}{\cos\frac\pi{12}}=
\sqrt2\,\big(\sqrt3-1\big)\cos\frac\pi2\Big(a+\frac16\Big).
\ee
When $a=\frac13$, we have $R_3(\frac13)=\sqrt3-1$, in agreement with
the first ratio in \eqref{BRR}.

\medskip
The ratio $R_3(a)$ in \eqref{BF3} is a remarkable result.
It combines one-parameter evaluations of a non-trivial $_2F_1$ function at two
different arguments of the form $\alpha\pm\beta$, which can be expressed as cubes of trigonometric functions, to yield a periodic cosine function.
The plots of the functions involved in \eqref{BF3} and their ratio $R_3$
as functions of the parameter $a$ are shown in Fig.~\ref{RRP}.

\section{Summary and outlook}

In the present communication, we have revised and extended the domain of validity for the celebrated hypergeometric transformation of Ramanujan, Berndt, Bhargava, and Garvan \eqref{RRR}.
Evaluations of the Gauss hypergeometric functions involved in \eqref{RRR} were performed at the endpoints of the established validity domain.

We have shown how the evaluations of Berndt and Chan related to \eqref{RRR} are connected with the singular value theory of complete elliptic integrals of the first kind.

For the more general $_2F_1$ evaluations of Ebisu and Beukers--Forsg{\aa}rd, which
depend on a free parameter $a$, the purely algebraic oscillatory ratios were calculated and visualized.

It is conceivable that similar one-parameter evaluations can be obtained directly from the classical hypergeometric transformation formulas of Kummer \cite{Kummer} and Goursat \cite{G1881}.

\section*{Acknowledgements}
Financial support through the research project ``Analytical and numerical approaches to analysis of cooperative behavior in complex systems'' by the National Academy of Sciences of Ukraine under the State Registration Number 0123U100238 is gratefully acknowledged. The author acknowledges the helpful assistance of Gemini by Google in final editing.\\
The author is deeply grateful to all warriors of the Ukrainian Armed Forces, living and fallen, for making it possible to continue the research work.

\providecommand{\href}[2]{#2}

\end{document}